\documentclass{scrartcl}
\usepackage{amsfonts}
\usepackage{amsmath}
\usepackage{amssymb}
\usepackage{pgf}
\usepackage{hyperref}
\usepackage{color}
\newcommand{\Idx}{\mathcal{I}}
\newcommand{\Jdx}{\mathcal{J}}
\newcommand{\chil}{\mathop{\operatorname{chil}}\nolimits}
\newcommand{\ctI}{\mathcal{T}_\Idx}
\newcommand{\ctJ}{\mathcal{T}_\Jdx}
\newcommand{\ctIJ}{\mathcal{T}_{\Idx\times\Jdx}}
\newcommand{\lfaIJ}{\mathcal{L}^+_{\Idx\times\Jdx}}
\newcommand{\lfiIJ}{\mathcal{L}^-_{\Idx\times\Jdx}}
\newcommand{\diam}{\mathop{\operatorname{diam}}\nolimits}
\newcommand{\dist}{\mathop{\operatorname{dist}}\nolimits}
\newcommand{\bbbn}{\mathbb{N}}
\newcommand{\bbbr}{\mathbb{R}}
\newtheorem{definition}{Definition}
\begin{document}
\title{Fast large-scale boundary element algorithms}
\author{Steffen B\"orm}
\maketitle              
\begin{abstract}
Boundary element methods (BEM) reduce a partial differential equation in
a domain to an integral equation on the domain's boundary.
They are particularly attractive for solving problems on unbounded
domains, but handling the dense matrices corresponding to the
integral operators requires efficient algorithms.

This article describes two approaches that allow us to solve
boundary element equations on surface meshes consisting of
several millions of triangles while preserving the optimal convergence
rates of the Galerkin discretization.
\end{abstract}

\section{Introduction}

We consider Laplace's equation
\begin{subequations}
\begin{align}\label{eq:laplace}
  \Delta u(x) &= 0 &
  &\text{ for all } x\in\Omega,
\end{align}
where $\Omega\subseteq\bbbr^3$ is a non-empty domain with a sufficiently
smooth boundary $\partial\Omega$.
If we add the boundary condition
\begin{align}\label{eq:dirichlet}
  u(x) &= f(x) &
  &\text{ for all } x\in\partial\Omega,
\end{align}
with a suitable function $f$, we obtain the \emph{Dirichlet problem}.
With the boundary condition
\begin{align}\label{eq:neumann}
  \frac{\partial u}{\partial n}(x) &= f(x) &
  &\text{ for all } x\in\partial\Omega,
\end{align}
\end{subequations}
where $n$ denotes the outward-pointing unit normal vector for the
domain $\Omega$, we arrive at the \emph{Neumann problem}.
We can solve these problems by using \emph{Green's representation
formula} (cf., e.g., \cite[Theorem~2.2.2]{HA92}) given by
\begin{align}\label{eq:green}
  u(x)
  &= \int_{\partial\Omega} g(x,y) \frac{\partial u}{\partial n}(y) \,dy
     - \int_{\partial\Omega} \frac{\partial g}{\partial n_y}(x,y) u(y) \,dy &
  &\text{ for all } x\in\Omega,
\end{align}
where
\begin{align*}
  g(x,y) &= \begin{cases}
    \frac{1}{4\pi \|x-y\|} &\text{ if } x\neq y,\\
    0 &\text{ otherwise}
  \end{cases} &
  &\text{ for all } x,y\in\bbbr^3
\end{align*}
denotes the free-space Green's function for the Laplace operator.
If we know Dirichlet and Neumann boundary conditions, we can use
(\ref{eq:green}) to compute the solution $u$ in any point of the
domain $\Omega$.

For boundary values $x\in\partial\Omega$, Green's formula takes the
form
\begin{align}
  \frac{1}{2} u(x)
  &= \int_{\partial\Omega} g(x,y) \frac{\partial u}{\partial n}(y) \,dy
   - \int_{\partial\Omega} \frac{\partial g}{\partial n_y}(x,y) u(y) \,dy
   \label{eq:dirichlet_to_neumann}\\
  &\qquad\text{ for almost all } x\in\partial\Omega,\notag
\end{align}
at least in the distributional sense \cite[eq. (3.92)]{SASC11},
and this boundary integral equation can be used to solve the Dirichlet
problem: we can solve the integral equation to obtain the Neumann
values $\frac{\partial u}{\partial n}$ from the Dirichlet values and
then use (\ref{eq:green}) to find the solution $u$ in all of $\Omega$.

In order to solve the Neumann problem, we take the normal derivative
of (\ref{eq:dirichlet_to_neumann}) and arrive at the boundary integral
equation
\begin{align}
  \frac{1}{2} \frac{\partial u}{\partial n}(x)
  &= \frac{\partial}{\partial n_x}
     \int_{\partial\Omega} g(x,y) \frac{\partial u}{\partial n}(y) \,dy
   - \frac{\partial}{\partial n_x}
     \int_{\partial\Omega} \frac{\partial g}{\partial n_y}(x,y) u(y) \,dy
   \label{eq:neumann_to_dirichlet}\\
  &\qquad\text{ for almost all } x\in\partial\Omega,\notag
\end{align}
again in the distributional sense \cite[eq. (3.92)]{SASC11},
that can be used to find Dirichlet values matching given Neumann
values, so we can follow the same approach as for the Dirichlet
problem.
The Dirichlet values are only determined up to a constant by
the Neumann values, and factoring out the constants leads to
unique solutions \cite[Theorem~3.5.3]{SASC11}.

In this paper, we will concentrate on solving the boundary integral
equations (\ref{eq:dirichlet_to_neumann}) and (\ref{eq:neumann_to_dirichlet})
efficiently.
We employ a Galerkin scheme using a finite-dimensional space
$V_h$ spanned by basis functions $(\psi_i)_{i\in\Idx}$ for the Neumann
values and another finite-dimensional space $U_h$ spanned by basis
functions $(\varphi_j)_{j\in\Jdx}$ for the Dirichlet values.
The discretization turns the boundary integral operators into
matrices, i.e., the matrix $G\in\bbbr^{\Idx\times\Idx}$ corresponding
to the \emph{single-layer operator} given by
\begin{subequations}
\begin{align}
  g_{ij} &= \int_{\partial\Omega} \psi_i(x)
           \int_{\partial\Omega} g(x,y) \psi_j(y) \,dy\,dx &
  &\text{ for all } i,j\in\Idx,\label{eq:slp}
\intertext{the matrix $K\in\bbbr^{\Idx\times\Jdx}$ corresponding to
the \emph{double-layer operator} given by}
  k_{ij} &= \int_{\partial\Omega} \psi_i(x)
           \int_{\partial\Omega} \frac{\partial g}{\partial n_y}(x,y)
           \varphi_j(y) \,dy \,dx &
  &\text{ for all } i\in\Idx,\ j\in\Jdx,\label{eq:dlp}
\intertext{and the matrix $W\in\bbbr^{\Jdx\times\Jdx}$ corresponding
to the \emph{hypersingular operator} given by}
  w_{ij} &= -\int_{\partial\Omega} \varphi_i(x) \frac{\partial}{\partial n_x}
           \int_{\partial\Omega} \frac{\partial g}{\partial n_y}(x,y)
           \varphi_j(y) \,dy \,dx &
  &\text{ for all } i,j\in\Jdx.
\end{align}
\end{subequations}
In order to set up this last matrix, we use an alternative
representation \cite[Corollary~3.3.24]{SASC11} based on the
single-layer operator.
Together with the mixed mass matrix $M\in\bbbr^{\Idx\times\Jdx}$ given by
\begin{align*}
  m_{ij} &= \int_{\partial\Omega} \psi_i(x) \varphi_j(y) &
  &\text{ for all } i\in\Idx,\ j\in\Jdx,
\end{align*}
we obtain the linear system
\begin{equation}\label{eq:matrix_dtn}
  G x = \left( \frac{1}{2} M + K \right) b
\end{equation}
for the Dirichlet-to-Neumann problem (\ref{eq:dirichlet_to_neumann}),
where $x\in\bbbr^\Idx$ contains the coefficients of the Neumann values
and $b\in\bbbr^\Jdx$ those of the given Dirichlet values, and the
system
\begin{equation}\label{eq:matrix_ntd}
  W x = \left( \frac{1}{2} M^* - K^* \right) b
\end{equation}
for the Neumann-to-Dirichlet problem (\ref{eq:neumann_to_dirichlet}),
where $x\in\bbbr^\Jdx$ now contains the coefficients of the Dirichlet
values and $b\in\bbbr^\Idx$ those of the given Neumann values.
$M^*$ and $K^*$ denote the transposed matrices of $M$ and $K$,
respectively.

Considered from the point of view of numerical mathematics, these
linear systems pose three challenges:
we have to evaluate singular integrals in order to compute diagonal
and near-diagonal entries, the resulting matrices are typically
dense and large, and the condition number grows with the matrix
dimension, which leads to a deteriorating performance of standard
Krylov solvers.
The first challenge can be met by using suitable quadrature
schemes like the Sauter-Schwab-Erichsen technique \cite{SA96,ERSA98,SASC11}.
For the third challenge, various preconditioning algorithms have
been proposed \cite{STWE98,LAPURE03}, among which we choose
$\mathcal{H}$-matrix coarsening and factorization \cite{GR04,BE05b}.

This leaves us with the second challenge, i.e., the efficient
representation of the dense matrices $G$, $K$, and $W$.
This task is usually tackled by employing a data-sparse approximation,
i.e., by finding a sufficiently accurate approximation of the matrices
that requires only a small amount of data.
One possibility to construct such an approximation is to use
wavelet basis functions and neglect small matrix entries in order
to obtain a sparse matrix \cite{PESC97,TA02,CODADE01,HASC06}.
Since the construction of suitable wavelet spaces on general surfaces
is complicated \cite{DASC99}, we will not focus on this approach.

Instead, we consider low-rank approximation techniques that directly
provide an approximation of the matrices for standard finite element
basis functions.
These techniques broadly fall into three categories:
\emph{analytic} methods approximate the kernel function $g$ locally
by sums of tensor products, and discretization of these sums
leads to low-rank approximations of matrix blocks.
The most prominent analytic methods are the fast multipole expansion
\cite{RO85,GRRO87,GRRO97} and interpolation \cite{GI01,BOGR02}.
\emph{Algebraic} methods, on the other hand, directly approximate
the matrix blocks, e.g., by computing a rank-revealing factorization
like an \emph{adaptive cross approximation} \cite{BE00a,BERJ01,BEKUVE15}.
The convergence and robustness of analytic methods can be proven
rigorously, but they frequently require a larger than necessary
amount of storage.
Algebraic methods reach close to optimal compression, but involve
a heuristic pivoting strategy that may fail in certain cases
\cite[Example~2.2]{BOGR04}.
\emph{Hybrid} methods combine analytic and algebraic techniques
in order to obtain the advantages of both without the respective
disadvantages.
In this paper, we will focus on the \emph{hybrid cross approximation}
\cite{BOGR04} and the \emph{Green cross approximation} \cite{BOCH14}
that both combine an analytic approximation with an algebraic
improvement in order to obtain fast and reliable algorithms.

\section{\texorpdfstring{$\mathcal{H}^2$-matrices}{H2-matrices}}

Both approximation schemes in this paper's focus lead to
\emph{$\mathcal{H}^2$-matrices} \cite{HAKHSA00,BO10}, a special
case of \emph{hierarchical matrices} \cite{HA99,GRHA02}.
A given matrix $G\in\bbbr^{\Idx\times\Jdx}$ with general finite
row and column index sets $\Idx$ and $\Jdx$ is split into
submatrices that are approximated by factorized low-rank
representations.

The submatrices are constructed hierarchically in order to
make the matrix accessible for elegant and efficient recursive
algorithms.
The first step is to split the index sets into a tree structure
of subsets.

%
%
\begin{definition}[Cluster tree]
Let $\Idx$ be a finite non-empty index set.
A tree $\mathcal{T}$ is called a \emph{cluster tree} for $\Idx$
if the following conditions hold:
\begin{itemize}
  \item Every node of the tree is a subset of $\Idx$.
  \item The root is $\Idx$.
  \item If a node has children, it is the union of these children:
    \begin{align*}
      t &= \bigcup_{t'\in\chil(t)} t' &
      &\text{ for all } t\in\mathcal{T} \text{ with } \chil(t)\neq\emptyset.
    \end{align*}
  \item The children of a node $t\in\mathcal{T}$ are disjoint:
    \begin{align*}
      t_1 \cap t_2 \neq \emptyset &\Rightarrow t_1=t_2 &
      &\text{ for all } t\in\mathcal{T},\ t_1,t_2\in\chil(t).
    \end{align*}
\end{itemize}
Nodes of a cluster tree are called \emph{clusters}.
\end{definition}

Cluster trees can be constructed by recursively splitting
index sets, e.g., ensuring that ``geometrically close'' indices
are contained in the same cluster \cite[Section~5.4]{HA15}.
We assume that cluster trees $\ctI$ and $\ctJ$ for the index
sets $\Idx$ and $\Jdx$ are given.

Using the cluster trees, the index set $\Idx\times\Jdx$ corresponding
to the matrix $G$ can now be split into a tree structure.

%
%
\begin{definition}[Block tree]
A tree $\mathcal{T}$ is called a \emph{block tree} for the row
index set $\Idx$ and the column index set $\Jdx$ if the following
conditions hold:
\begin{itemize}
  \item Every node of the tree is a subset
    $t\times s\subseteq\Idx\times\Jdx$ with $t\in\ctI$ and $s\in\ctJ$.
  \item The root is $\Idx\times\Jdx$.
  \item If a node has children, the children are given as follows:
    \begin{align*}
      \chil(t\times s) &= \begin{cases}
        \{ t\times s'\ :\ s'\in\chil(s) \}
        &\text{ if } \chil(t)=\emptyset,\\
        \{ t'\times s\ :\ t'\in\chil(t) \}
        &\text{ if } \chil(s)=\emptyset,\\
        \{ t'\times s'\ :\ t'\in\chil(t),\ s'\in\chil(s) \}
        &\text{ otherwise}
      \end{cases}
    \end{align*}
    for all $t\times s\in\mathcal{T}$.
\end{itemize}
Nodes of a block tree are called \emph{blocks}.
\end{definition}

Block trees are usually constructed recursively using an
\emph{admissibility condition} matching the intended approximation
scheme:
we start with the root $\Idx\times\Jdx$ and recursively subdivide
blocks.
Once the admissibility condition indicates for a block $t\times s$
that we can approximate the submatrix $G|_{t\times s}$, we stop
subdividing and call the block $t\times s$ a \emph{farfield block}.
We also stop if $t$ and $s$ have no children, then $t\times s$
is a \emph{nearfield block}.
If we ensure that leaf clusters contain only a small number of
indices, nearfield blocks are small and we can afford to store them
without compression.
The key to the efficiency of hierarchical matrices is the data-sparse
representation of the farfield blocks.

We assume that a block tree $\ctIJ$ is given and that its farfield
leaves are collected in a set $\lfaIJ$, while the nearfield leaves
are in a set $\lfiIJ$.

For a hierarchical matrix \cite{HA99,GRHA02}, we simply assume that
farfield blocks have low rank $k\in\bbbn$ and can therefore be stored
efficiently in factorized form $G|_{t\times s} \approx A_{ts} B_{ts}^*$ with
$A_{ts}\in\bbbr^{t\times k}$, $B_{ts}\in\bbbr^{s\times k}$.
Here we use $\bbbr^{t\times k}$ and $\bbbr^{s\times k}$ as abbreviations
for $\bbbr^{t\times[1:k]}$ and $\bbbr^{s\times[1:k]}$.

The more efficient $\mathcal{H}^2$-matrices \cite{HAKHSA00,BO10}
use a different factorized representation closely related to
fast multipole methods \cite{RO85}:
each cluster is associated with a low-dimensional subspace, e.g., a
space spanned by polynomials or multipole functions, and the range of
a matrix block $G|_{t\times s}$ has to be contained in the subspace
for the row cluster $t$, while the range of the adjoint block
$G|_{t\times s}^*$ has to be contained in the subspace for the
column cluster $s$.
This property is expressed by the equation (\ref{eq:vsw}) below.

%
%
\begin{definition}[Cluster basis]
Let $k\in\bbbn$.
A family $V=(V_t)_{t\in\ctI}$ of matrices is called a \emph{cluster basis}
for $\ctI$ with (maximal) rank $k$ if the following conditions hold:
\begin{itemize}
  \item We have $V_t\in\bbbr^{t\times k}$ for all $t\in\ctI$.
  \item For all $t\in\ctI$ with $\chil(t)\neq\emptyset$, there are
    \emph{transfer matrices} $E_{t'}\in\bbbr^{k\times k}$ for all
    $t'\in\chil(t)$ such that
    \begin{align}\label{eq:transfer}
      V_t|_{t'\times k} &= V_{t'} E_{t'} &
      &\text{ for all } t'\in\chil(t).
    \end{align}
\end{itemize}
\end{definition}

An important property of cluster bases is that they can be stored
efficiently:
under standard assumptions, only $\mathcal{O}(n k)$ units of
storage are requires, where $n=|\Idx|$ is the cardinality of
the index set $\Idx$.

%
%
\begin{definition}[$\mathcal{H}^2$-matrix]
Let $V=(V_t)_{t\in\ctI}$ and $W=(W_s)_{s\in\ctJ}$ be cluster bases for
$\ctI$ and $\ctJ$, respectively.
A matrix $G\in\bbbr^{\Idx\times\Jdx}$ is called an
\emph{$\mathcal{H}^2$-matrix} with row basis $(V_t)_{t\in\ctI}$
and column basis $(W_s)_{s\in\ctJ}$ if for every admissible
leaf $t\times s\in\lfaIJ$ there is a \emph{coupling matrix}
$S_{ts}\in\bbbr^{k\times k}$ such that
\begin{equation}\label{eq:vsw}
  G|_{t\times s} = V_t S_{ts} W_s^*.
\end{equation}
\end{definition}

Under standard assumptions, we can represent an $\mathcal{H}^2$-matrix
by $\mathcal{O}(nk+mk)$ coefficients, where $n=|\Idx|$ and $m=|\Jdx|$
are the cardinalities of the index sets.
The matrix-vector multiplication $x\mapsto G x$ can be performed in
$\mathcal{O}(n k + m k)$ operations for $\mathcal{H}^2$-matrices, and
there are a number of other important operations that also only have
linear complexity with respect to $n$ and $m$, cf. \cite{BO10}.

\section{Hybrid cross approximation}

In order to construct an $\mathcal{H}^2$-matrix approximation of
the matrices $V$, $K$, and $W$ required for the boundary element
method, we first consider the \emph{hybrid cross approximation}
(HCA) method \cite{BOGR04} originally developed for hierarchical
matrices.

We associate each cluster $t\in\ctI$ with an axis-parallel
\emph{bounding box} $B_t\subseteq\bbbr^3$ such that the supports
of all basis functions associated with indices in $t$ are
contained in $B_t$.
In order to ensure that the kernel function $g$ is sufficiently
smooth for a polynomial approximation, we introduce the
admissibility condition
\begin{equation}\label{eq:admissibility}
  \max\{ \diam(B_t), \diam(B_s) \} \leq 2 \eta \dist(B_t, B_s),
\end{equation}
where $\diam(B_t)$ and $\diam(B_s)$ denote the Euclidean diameters
of the bounding boxes $B_t$ and $B_s$, while $\dist(B_t, B_s)$
denotes their Euclidean distance.
$\eta$ is a parameter that controls the storage complexity and
the accuracy of the approximation.
In our experiments, the choice $\eta=1$ leads to reasonable
results.

If two clusters $t\in\ctI$, $s\in\ctJ$ satisfy this condition,
the restriction $g|_{B_t\times B_s}$ can be approximated by
polynomials.
We use $m$-th order tensor Chebyshev interpolation and denote the
interpolation points for $B_t$ and $B_s$ by $(\xi_{t,\nu})_{\nu=1}^k$ and
$(\xi_{s,\mu})_{\mu=1}^k$, where $k=m^3$.
The corresponding Lagrange polynomials are denoted by
$(\mathcal{L}_{t,\nu})_{\nu=1}^k$ and $(\mathcal{L}_{s,\mu})_{\mu=1}^k$,
and the tensor interpolation polynomial is given by
\begin{align*}
  \tilde g_\text{int}(x,y) &= \sum_{\nu=1}^k \sum_{\mu=1}^k
    \mathcal{L}_{t,\nu}(x)
    g(\xi_{t,\nu}, \xi_{s,\mu})
    \mathcal{L}_{s,\mu}(y) &
  &\text{ for all } x\in B_t,\ y\in B_s.
\end{align*}
It is possible to prove
\begin{align*}
  \|g - \tilde g_\text{int}\|_{\infty,B_t\times B_s}
  &\lesssim \frac{q^m}{\diam(B_t)^{1/2} \diam(B_s)^{1/2}}\\
  &\text{ for all } m\in\bbbn
   \text{ and all } t\in\ctI, s\in\ctJ
     \text{ satisfying (\ref{eq:admissibility})}.
\end{align*}
The rate $q$ of convergence depends only on the parameter $\eta$,
cf. \cite[Chapter~4]{BO10}.

Unfortunately, the rank $k$ of this approximation is too large:
potential theory suggests that a rank $k\sim m^2$ should be
sufficient for an $m$-th order approximation.

We reduce the rank by combining the interpolation with an
algebraic procedure, in this case \emph{adaptive cross approximation}
\cite{BE00a,BERJ01}:
we introduce the matrix $S\in\bbbr^{k\times k}$ by
\begin{align*}
  s_{\nu\mu} &= g(\xi_{t,\nu}, \xi_{s,\mu}) &
  &\text{ for all } \nu,\mu\in[1:k]
\end{align*}
and use a rank-revealing pivoted LU factorization to obtain an
approximation of the form
\begin{equation*}
  S \approx S|_{[1:k]\times\sigma} C
            S|_{\tau\times[1:k]},
\end{equation*}
where $C:=(S|_{\tau\times\sigma})^{-1}$ and $\tau,\sigma\subseteq[1:k]$
denote the first $\tilde k\leq k$ row and column pivots, respectively.
Due to the properties of the kernel function $g$, the rank
$\tilde k$ can be expected to be significantly smaller than $k$.

Given the pivot sets $\tau$ and $\sigma$, we can now ``take back''
the interpolation:
\begin{align*}
  g(x,y) &\approx \sum_{\nu=1}^k \sum_{\mu=1}^k
     \mathcal{L}_{t,\nu}(x) s_{\nu\mu} \mathcal{L}_{s,\mu}(y)\\
  &\approx \sum_{\nu=1}^k \sum_{\mu=1}^k
     \mathcal{L}_{t,\nu}(x) (S|_{[1:k]\times\sigma} C
     S|_{\tau\times[1:k]})_{\nu\mu} \mathcal{L}_{s,\mu}(y)\\
  &= \sum_{\nu=1}^k \sum_{\mu=1}^k
     \sum_{\lambda\in\tau} \sum_{\kappa\in\sigma}
     \mathcal{L}_{t,\nu}(x) s_{\nu\kappa} c_{\kappa\lambda}
     s_{\lambda\mu} \mathcal{L}_{s,\mu}(y)\\
  &= \sum_{\lambda\in\tau} \sum_{\kappa\in\sigma}
     \underbrace{\sum_{\nu=1}^k \mathcal{L}_{t,\nu}(x)
                    g(\xi_{t,\nu}, \xi_{s,\kappa})
                }_{\approx g(x,\xi_{s,\kappa})}
     c_{\kappa\lambda}
     \underbrace{\sum_{\mu=1}^k g(\xi_{t,\lambda}, \xi_{s,\mu})
                    \mathcal{L}_{s,\mu}(y)
                }_{\approx g(\xi_{t,\lambda},y)}\\
  &\approx \sum_{\lambda\in\tau} \sum_{\kappa\in\sigma}
     g(x,\xi_{s,\kappa})\ c_{\kappa\lambda}\ g(\xi_{t,\lambda},y).
     =: \tilde g_\text{hca}(x,y)
\end{align*}
By controlling the interpolation order and the accuracy of the
cross approximation, we can ensure that the approximation error
$\|g-\tilde g_\text{hca}\|_{\infty,B_t\times B_s}$ is below any given
tolerance \cite{BOGR04}.

In order to obtain an approximation of the submatrix
$G|_{t\times s}$, we replace $g$ by $\tilde g_\text{hca}$ in (\ref{eq:slp})
to find
\begin{align}
  g_{ij} &= \int_{\partial\Omega} \psi_i(x)
           \int_{\partial\Omega} g(x,y) \psi_j(y) \,dy\,dx
   \notag\\
  &\approx \int_{\partial\Omega} \psi_i(x)
           \int_{\partial\Omega} \tilde g_\text{hca}(x,y) \psi_j(y) \,dy\,dx
   \notag\\
  &= \sum_{\lambda\in\tau} \sum_{\kappa\in\sigma}
     \underbrace{\int_{\partial\Omega} \psi_i(x) g(x,\xi_{s,\kappa})
                }_{=: a_{i\kappa}}
     c_{\lambda\kappa}
     \underbrace{\int_{\partial\Omega} \psi_j(y) g(\xi_{t,\lambda},y)
                }_{=: b_{j\lambda}}\label{eq:hca_discrete}\\
  &= (A C B^*)_{ij}
    \qquad\text{ for all } i\in t,\ j\in s\notag
\end{align}
with matrices $A\in\bbbr^{t\times\sigma}$ and $B\in\bbbr^{s\times\tau}$.
Due to $|\tau|=|\sigma|=\tilde k\leq k$, the matrix $A C B^*$ is an
improved low-rank approximation of the submatrix $G|_{t\times s}$.

This procedure alone only yields a hierarchical matrix, not the
desired $\mathcal{H}^2$-matrix.
In order to reduce the storage requirements, we apply
\emph{hierarchical compression} \cite{BO07a}: the hybrid cross
approximation technique already yields low-rank approximations
of individual blocks, but these approximations do not share
common row or column cluster bases.
The hierarchical compression algorithm recursively merges independent
submatrices into larger $\mathcal{H}^2$-matrices until the entire
matrix is in the required form.

The parallelization of this algorithm is fairly straightforward:
we can set up the matrices corresponding to the leaves of the block
tree in parallel, and the merge operations for submatrices can also
be performed in parallel once their children are available.
In order to handle the dependencies between children and their
parents, both clusters and blocks, a task-based programming model
is particularly useful for this algorithm.

\section{Green cross approximation}

While hybrid cross approximation requires a subsequent compression
step to obtain an $\mathcal{H}^2$-matrix, we now consider a related
technique that \emph{directly} yields an $\mathcal{H}^2$-matrix
approximation of the matrix $G$.

\begin{figure}
  \pgfdeclareimage[width=7cm]{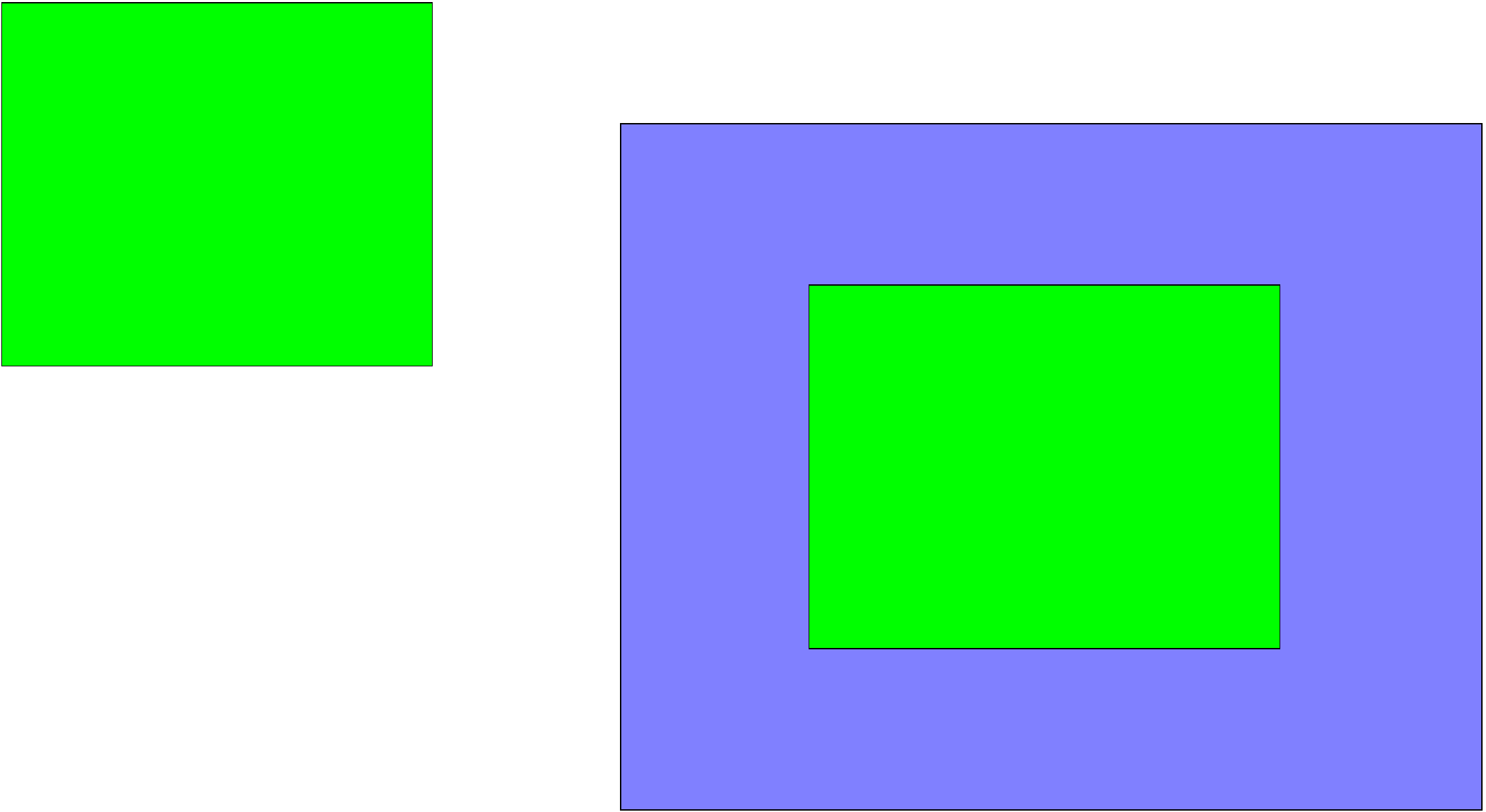}{omega_t}

  \begin{center}
    \begin{pgfpicture}{0cm}{0cm}{7cm}{4cm}
      \pgfuseimage{omega_t}%
      \pgfputat{\pgfxy(-2.1,1.5)}{\pgfbox[center,base]{$B_t$}}%
      \pgfputat{\pgfxy(-6.0,2.8)}{\pgfbox[center,base]{$B_s$}}%
      \pgfputat{\pgfxy(-2.1,3.4)}{\pgfbox[center,base]{$\partial\omega_t$}}%
    \end{pgfpicture}
  \end{center}

  \caption{Possible choice of the auxiliary bounding box $\omega_t$
    corresponding to $B_t$ and $B_s$}
  \label{fi:omega_t}
\end{figure}

We once more consider an admissible pair $t\in\ctI$, $s\in\ctJ$
of clusters with bounding boxes $B_t,B_s\subseteq\bbbr^3$.
We construct an auxiliary axis-parallel box $\omega_t\subseteq\bbbr^3$
such that
\begin{align}\label{eq:omega_t}
  B_t &\subseteq \omega_t, &
  \diam(\omega_t) &\lesssim \dist(B_t,\partial\omega_t), &
  \diam(\omega_t) &\lesssim \dist(B_s,\omega_t),
\end{align}
cf. Figure~\ref{fi:omega_t} for an illustration.
Let $y\in B_s$.
The third assumption in (\ref{eq:omega_t}) implies $y\not\in\omega_t$,
therefore the function $u(x) = g(x,y)$ is harmonic in $\omega_t$
This property allows us to apply Green's equation (\ref{eq:green})
to the domain $\omega_t$ to obtain
\begin{equation*}
  g(x,y)
  = \int_{\partial\omega_t} g(x,z) \frac{\partial g}{\partial n_z}(z,y) \,dz
  - \int_{\partial\omega_t} \frac{\partial g}{\partial n_z}(x,z) g(z,y) \,dz.
\end{equation*}
We observe that the variables $x$ and $y$ are separated in both
integrands.
Due to the second and third assumptions in (\ref{eq:omega_t}), the
integrands are smooth for $x\in B_t$ and $y\in B_s$, and we can approximate
the integrals by a quadrature rule, e.g., a composite Gauss rule, with weights
$(w_\nu)_{\nu=1}^k$ and quadrature points $(z_\nu)_{\nu=1}^k$ on
the boundary $\partial\omega_t$ in order to get
\begin{equation*}
  g(x,y)
  \approx \sum_{\nu=1}^k w_\nu
     g(x,z_\nu) \frac{\partial g}{\partial n_z}(z_\nu,y)
  - \sum_{\nu=1}^k w_\nu
     \frac{\partial g}{\partial n_z}(x,z_\nu) g(z_\nu,y)
  =: \tilde g_\text{grn}(x,y)
\end{equation*}
for all $x\in B_t$ and all $y\in B_s$.
The function $\tilde g_\text{grn}$ is again a sum of tensor
products, and, as in the case of the hybrid cross approximation,
its discretization gives rise to a low-rank approximation of
the submatrix $G|_{t\times s}$, since we have
\begin{align*}
  g_{ij}
  &= \int_{\partial\Omega} \psi_i(x)
     \int_{\partial\Omega} g(x,y) \psi_j(y) \,dy \,dx\\
  &\approx \int_{\partial\Omega} \psi_i(x)
     \int_{\partial\Omega} \tilde g_\text{grn}(x,y) \psi_j(y) \,dy \,dx\\
  &= \sum_{\nu=1}^k
     \underbrace{w_\nu^{1/2} \int_{\partial\Omega} \psi_i(x) g(x,z_\nu) \, dx
                }_{=:a_{i\nu}}
     \underbrace{w_\nu^{1/2} \int_{\partial\Omega} \psi_j(y)
          \frac{\partial g}{\partial n_z}(z_\nu,y) \,dy
                }_{=:b_{j\nu}}\\
  &- \sum_{\nu=1}^k
     \underbrace{w_\nu^{1/2} \int_{\partial\Omega} \psi_i(x)
          \frac{\partial g}{\partial n_z}(x,z_\nu) \,dx
                }_{=:c_{i\nu}}
     \underbrace{w_\nu^{1/2} \int_{\partial\Omega} \psi_j(x) g(z_\nu,y) \,dy
                }_{=:d_{j\nu}}\\
  &= (A B^* - C D^*)_{ij}
   \qquad\text{ for all } i\in t,\ j\in s
\end{align*}
with $A,C\in\bbbr^{t\times k}$ and $B,D\in\bbbr^{s\times k}$, therefore
\begin{equation*}
  G|_{t\times s}
  \approx \begin{pmatrix}
    A & C
  \end{pmatrix}
  \begin{pmatrix}
    B^*\\ D^*
  \end{pmatrix}.
\end{equation*}
Unfortunately, the rank of this approximation is quite high, higher
than, e.g., for the hybrid cross approximation.
Once again, we can use an algebraic technique, i.e., the adaptive
cross approximation, to improve the analytically-motivated initial
approximation $\tilde g_\text{grn}$:
we define
\begin{equation*}
  L := \begin{pmatrix}
         A & C
       \end{pmatrix} \in \bbbr^{t\times (2k)}
\end{equation*}
and perform an adaptive cross approximation to find index sets
$\tau\subseteq t$ and $\sigma\subseteq [1:2k]$ of cardinality
$\tilde k\leq 2k$ such that
\begin{equation*}
  L \approx L|_{t\times\sigma} (L|_{\tau\times\sigma})^{-1} L|_{\tau\times [1:2k]}.
\end{equation*}
By applying this approximation and ``taking back'' the quadrature
approximation in the last step, we arrive at
\begin{equation*}
  G|_{t\times s}
  \approx L \begin{pmatrix}
              B^*\\ D^*
            \end{pmatrix}
  \approx L|_{t\times\sigma} (L|_{\tau\times\sigma})^{-1}
     L|_{\tau\times [1:2k]} \begin{pmatrix}
                          B^*\\ D^*
                        \end{pmatrix}
  \approx L|_{t\times\sigma} (L|_{\tau\times\sigma})^{-1}
     G|_{\tau\times s}.
\end{equation*}
We define $\hat t := \tau$ and $V_t := L|_{t\times\sigma}
(L|_{\tau\times\sigma})^{-1}\in\bbbr^{t\times\hat t}$ and obtain
\begin{equation}\label{eq:gca_onesided}
  G|_{t\times s} \approx V_t G|_{\hat t\times s}.
\end{equation}
This is a rank-$\tilde k$ approximation of the matrix block, and
experiments indicate that $\tilde k$ is frequently far smaller
than $2k$.
The equation (\ref{eq:gca_onesided}) can be interpreted as
``algebraic interpolation'': we (approximately) recover all entries
of the matrix $G|_{t\times s}$ from a few rows $G|_{\hat t\times s}$,
where the indices in $\hat t$ play the role of interpolation points
and the columns of $V_t$ the role of Lagrange polynomials.
We call this approach \emph{Green cross approximation} (GCA).

Concerning our goal of finding an $\mathcal{H}^2$-matrix, we
observe that $V_t$ and $\hat t$ depend only on $t$ and $\omega_t$,
but not on the cluster $s$, therefore $V_t$ it is a good candidate
for a cluster basis.

Applying the same procedure to the cluster $s$ instead of $t$,
we obtain $\hat s\subseteq s$ and $W_s\in\bbbr^{s\times\hat s}$ such
that
\begin{equation*}
  G|_{t\times s} \approx G|_{t\times\hat s} W_s^*,
\end{equation*}
and combining both approximations yields
\begin{equation*}
  G|_{t\times s} \approx V_t G|_{\hat t\times s}
               \approx V_t G|_{\hat t\times\hat s} W_s^*,
\end{equation*}
the required representation (\ref{eq:vsw}) for an $\mathcal{H}^2$-matrix.

In order to make $V=(V_t)_{t\in\ctI}$ and $W=(W_s)_{s\in\ctJ}$ proper
cluster bases, we have to ensure the nesting property
(\ref{eq:transfer}).
We can achieve this goal by slightly modifying our construction:
we assume that $t\in\ctI$ has two children $t_1,t_2\in\chil(t)$
and that the sets $\hat t_1\subseteq t_1$ and $\hat t_2\subseteq t_2$
and the matrices $V_{t_1}$ and $V_{t_2}$ have already been computed.
We have
\begin{equation*}
  G|_{t\times s}
  = \begin{pmatrix}
      G|_{t_1\times s}\\
      G|_{t_2\times s}
    \end{pmatrix}
  \approx \begin{pmatrix}
      V_{t_1} G|_{\hat t_1\times s}\\
      V_{t_2} G|_{\hat t_2\times s}
    \end{pmatrix}
  = \begin{pmatrix}
      V_{t_1} & \\
      & V_{t_2}
    \end{pmatrix}
    G|_{(\hat t_1\cup\hat t_2)\times s}.
\end{equation*}
We apply the cross approximation to $G|_{(\hat t_1\cup\hat t_2)\times s}$
instead of $G|_{t\times s}$ and obtain $\hat t\subseteq\hat t_1\cup\hat t_2$
and $\widehat V_t\in\bbbr^{(\hat t_1\cup\hat t_2)\times\hat t}$ such that
\begin{equation*}
  G|_{(\hat t_1\cup\hat t_2)\times s}
  \approx \widehat V_t G|_{\hat t\times s}
\end{equation*}
and therefore
\begin{equation*}
  G|_{t\times s}
  \approx \begin{pmatrix}
            V_{t_1} & \\
            & V_{t_2}
          \end{pmatrix} G|_{(\hat t_1\cup\hat t_2)\times s}
  \approx \begin{pmatrix}
            V_{t_1} & \\
            & V_{t_2}
          \end{pmatrix} \widehat V_t G|_{\hat t\times s},
\end{equation*}
so defining
\begin{equation*}
  V_t := \begin{pmatrix}
           V_{t_1} & \\
           & V_{t_2}
         \end{pmatrix} \widehat V_t
\end{equation*}
ensures (\ref{eq:transfer}) if we let
\begin{equation*}
  \begin{pmatrix}
    E_{t_1}\\ E_{t_2}
  \end{pmatrix} := \widehat V_t.
\end{equation*}
This modification not only ensures that $V=(V_t)_{t\in\ctI}$ is
a proper cluster basis, it also reduces the computational work
required for the cross approximation, since only the submatrices
$G|_{(\hat t_1\cup\hat t_2)\times s}$ have to be considered, and these
submatrices are significantly smaller than $G|_{t\times s}$.

The parallelization of this algorithm is straightforward:
we can compute index sets and matrices for all leaves of the
cluster tree in parallel.
Once this task has been completed, we can treat all clusters whose
children are leaves in parallel.
Next are all clusters whose children have already been treated.
Repeating this procedure until we reach the root of the tree
yields a simple and efficient parallel version of the algorithm.
Our implementation uses simple parallel for loops to iterate
through the children clusters and blocks prior to setting up
their parents.

\section{Numerical experiments}

Now that we have two compression algorithms at our disposal,
we have to investigate how well they perform in practice.
While we can prove for both algorithms that they can reach
any given accuracy (disregarding rounding errors), we have to
see which accuracies are necessary in order to preserve the
theoretical convergence rates of the Galerkin discretization.

Until now, we have only seen HCA and GCA applied to the
matrix $G$ corresponding to the single-layer operator.
For the other two operators, we simply take the appropriate
derivatives of $\tilde g_\text{hca}$ and $\tilde g_\text{grn}$
and use them as approximations of the kernel functions.

Given a boundary element mesh, we choose discontinuous piecewise
constant basis functions for the Neumann values and continuous
piecewise linear basis functions for the Dirichlet values.
For a meshwidth of $h\in\bbbr_{>0}$, we expect theoretical
convergence rates of $\mathcal{O}(h)$ for the Neumann values
in the $L^2$ norm, $\mathcal{O}(h^{3/2})$ for the Neumann values
in the $H^{-1/2}$ norm and the Dirichlet values in the $H^{1/2}$
norm, and $\mathcal{O}(h^2)$ for the Dirichlet values in the
$L^2$ norm.

In a first experiment, we approximate the unit sphere
$\{x\in\bbbr^3\ :\ x_1^2+x_2^2+x_3^2=1\}$ by a sequence of
triangular meshes constructed by splitting the eight sides
of a double pyramid $\{x\in\bbbr^3\ :\ |x_1|+|x_2|+|x_3|=1\}$
regularly into triangles and then projecting all vertices
to the unit sphere.
Since we expect the condition number of the matrices to be
in $\mathcal{O}(h^{-1})$ and want to preserve the theoretical
convergence rate of $\mathcal{O}(h^2)$, we aim for an
accuracy of $\mathcal{O}(h^3)$ for the matrix approximation.
For the sake of simplicity, we use a slightly higher accuracy:
if $h$ is halved, we reduce the error tolerance by a factor
of $10$ instead of just $8$.

Nearfield matrix entries are computed by
Sauter-Schwab-Erichsen quadrature \cite{SASC11}, and we have
to increase the order of the nearfield quadrature occasionally
to ensure the desired rate of convergence.

%
%
\begin{table}
  \caption{Parameters chosen for the unit sphere}
  \label{ta:sphere}

  \begin{equation*}
    \begin{array}{r|r|r|r|r|r|r|r}
      n & q_\text{near} & r_\text{leaf} & m
        & \epsilon_\text{aca} & \epsilon_\text{comp}
        & \epsilon_\text{slv} & \epsilon_\text{prc}\\
      \hline
      8\,192 & 4 & 25 & 5 & 1_{-5} & 1_{-5} & 1_{-6} & 1_{-2}\\
     18\,432 & 4 & 36 & 6 & 3_{-6} & 3_{-6} & 3_{-7} & 1_{-2}\\
     32\,768 & 4 & 36 & 6 & 1_{-6} & 1_{-6} & 1_{-7} & 5_{-3}\\
     73\,728 & 5 & 49 & 7 & 3_{-7} & 3_{-7} & 3_{-8} & 5_{-3}\\
    131\,072 & 5 & 49 & 7 & 1_{-7} & 1_{-7} & 1_{-8} & 2_{-3}\\
    294\,912 & 5 & 64 & 8 & 3_{-8} & 3_{-8} & 3_{-9} & 2_{-3}\\
    524\,288 & 5 & 64 & 8 & 1_{-8} & 1_{-8} & 1_{-9} & 1_{-3}\\
 1\,179\,648 & 6 & 81 & 9 & 3_{-9} & 3_{-9} & 3_{-10} & 1_{-3}\\
 2\,097\,152 & 6 & 81 & 9 & 1_{-9} & 1_{-9} & 1_{-10} & 5_{-4}\\
 4\,718\,592 & 6 & 100 & 10 & 3_{-10} & 3_{-10} & 3_{-11} & 5_{-4}\\
 8\,388\,608 & 6 & 100 & 10 & 1_{-10} & 1_{-10} & 1_{-11} & 2_{-4}
    \end{array}
  \end{equation*}
\end{table}

The parameters used for this experiment are summarized in
Table~\ref{ta:sphere}, where $n$ denotes the number of triangles,
$q_\text{near}$ the nearfield quadrature order, $r_\text{leaf}$ the
resolution of the cluster tree, i.e., the maximal size of leaf
clusters, $m$ the order of interpolation for HCA and the order
of quadrature for GCA, $\epsilon_\text{aca}$ the relative accuracy
for the adaptive cross approximation, $\epsilon_\text{comp}$ the
accuracy for the hierarchical compression used for HCA,
$\epsilon_\text{slv}$ the relative accuracy of the Krylov solver,
and $\epsilon_\text{prc}$ the relative accuracy of the
preconditioner constructed by coarsening \cite{GR04} and
$\mathcal{H}$-Cholesky decomposition \cite{HA15}.
We use $\eta=1$ for the admissibility parameter and construct
the boxes $\omega_t$ for the Green quadrature to ensure
$\dist(\partial\omega_t,B_t)=\delta_t$, with
$\delta_t=\max\{b_1-a_1,b_2-a_2,b_3-a_3\}$, where
$B_t=[a_1,b_1]\times[a_2,b_2]\times[a_3,b_3]$.
This choice is not strictly covered by the theorey in \cite{BOCH14},
but works well in practice.
The nearfield order was only increased if the convergence was
compromised.
$\epsilon_\text{aca}$, $\epsilon_\text{comp}$, and $\epsilon_\text{slv}$
where chosen in the expectation that an accuracy of $\mathcal{O}(h^3)$
would be required in order to keep up with the discretization error.
$\epsilon_\text{prc}$ was chosen in the expectation that an
accuracy of $\mathcal{O}(h)$ would be necessary to keep up with
the growth of the condition number of the linear system.

%
%
\begin{figure}
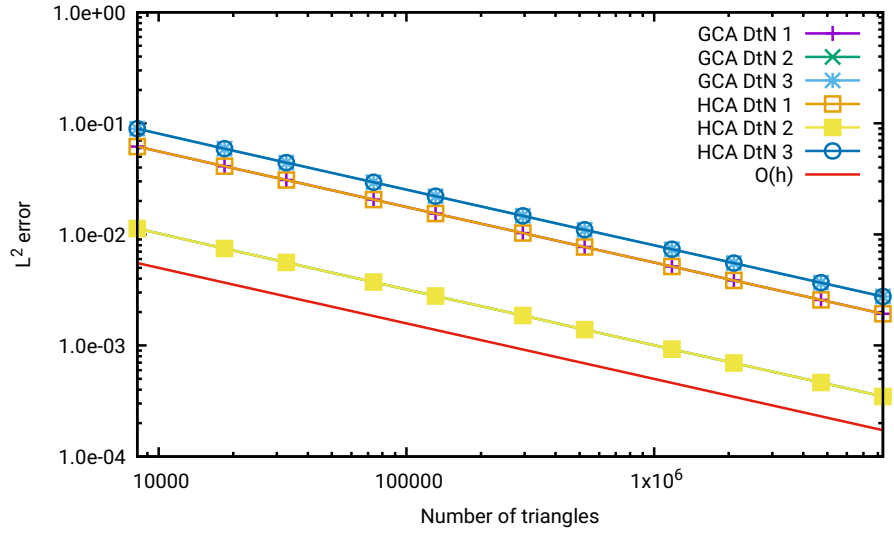

  \pgfdeclareimage[width=12cm]{l2_dtn}{gca_hca_cl_l2error_dtn}

  \begin{center}
    \pgfuseimage{l2_dtn}
  \end{center}

  \caption{$L^2$ error for the Dirichlet-to-Neumann problem}
  \label{fi:l2_dtn}
\end{figure}

Figure~\ref{fi:l2_dtn} shows the $L^2$-norm error for the
approximation of the Neumann data computed via the boundary
integral equation (\ref{eq:dirichlet_to_neumann}).
We use the functions $u_1(x)=x_1^2-x_3^2$, $u_2(x)=g(x,y_1)$,
and $u_3(x)=g(x,y_2)$ with $y_1=(1.2,1.2,1.2)$ and $y_2=(1.0,0.25,1.0)$
as test cases.
We can see that the optimal convergence rate of $\mathcal{O}(h)$
is preserved despite the matrix compression and nearfield quadrature.

%
%
\begin{figure}
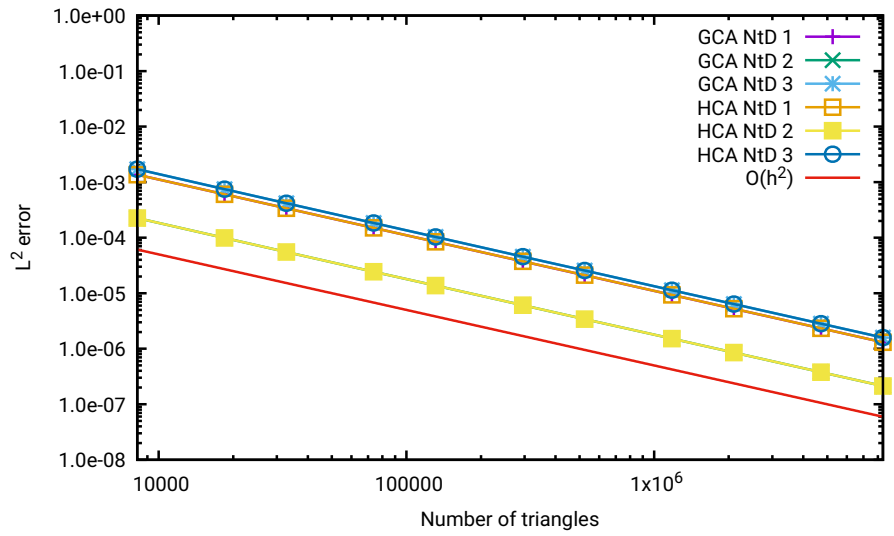

  \pgfdeclareimage[width=12cm]{l2_ntd}{gca_hca_cl_l2error_ntd}

  \begin{center}
    \pgfuseimage{l2_ntd}
  \end{center}

  \caption{$L^2$ error for the Neumann-to-Dirichlet problem}
  \label{fi:l2_ntd}
\end{figure}

Figure~\ref{fi:l2_ntd} shows the $L^2$-norm error for the
approximation of the Dirichlet data computed via the boundary
integral equation (\ref{eq:neumann_to_dirichlet}).
Also in this case, the optimal convergence rate of $\mathcal{O}(h^2)$
is preserved.

\begin{figure}
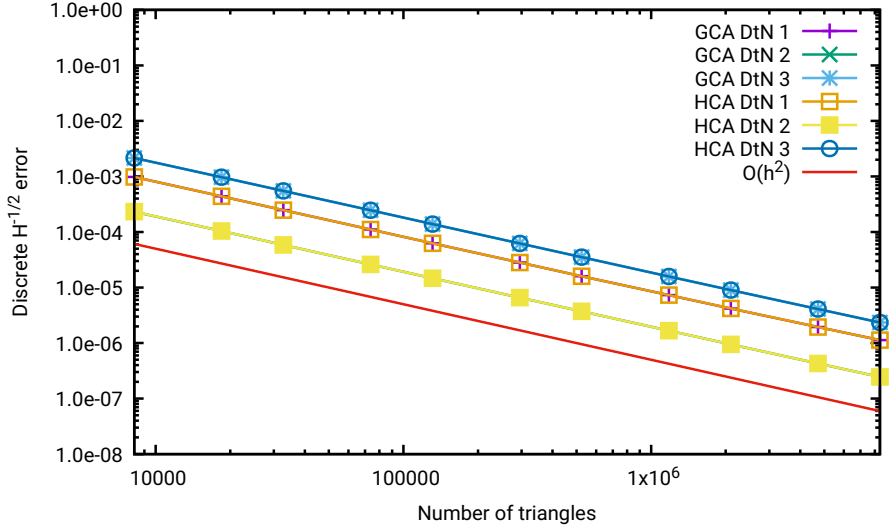

  \pgfdeclareimage[width=12cm]{energy_dtn}{gca_hca_cl_energyerror_dtn}

  \begin{center}
    \pgfuseimage{energy_dtn}
  \end{center}

  \caption{Discrete $H^{-1/2}$ error for the Dirichlet-to-Neumann problem}
  \label{fi:energy_dtn}
\end{figure}

Since computing the exact $H^{-1/2}$-norm error is complicated,
we rely on the $H^{-1/2}$-ellipticity of the single-layer operator:
we compute the $L^2$-projection of the Neumann values into the
discrete space, which is expected to converge at a rate of
$\mathcal{O}(h^{3/2})$ to the exact solution, and then compare it
to the Galerkin solution using the energy product corresponding
to the matrix $G$.
Figure~\ref{fi:energy_dtn} indicates that the discrete $H^{-1/2}$-norm
error even converges at a rate of $\mathcal{O}(h^2)$, therefore the
$\mathcal{O}(h^{3/2})$ convergence compared to the continuous solution
is also preserved.

%
%

Now that we have established that the compression algorithms do
not hurt the optimal convergence rates of the Galerkin discretization,
we can consider the corresponding complexity.
The runtimes have been measured on a system with two Intel\textregistered{}
Xeon\textregistered{} Platinum 8160 processors, each with 24 cores and
running at a base clock of $2.1$ GHz.
The implementation is based on the open-source H2Lib software package,
cf. \url{http://www.h2lib.org}.

%
%
\begin{figure}
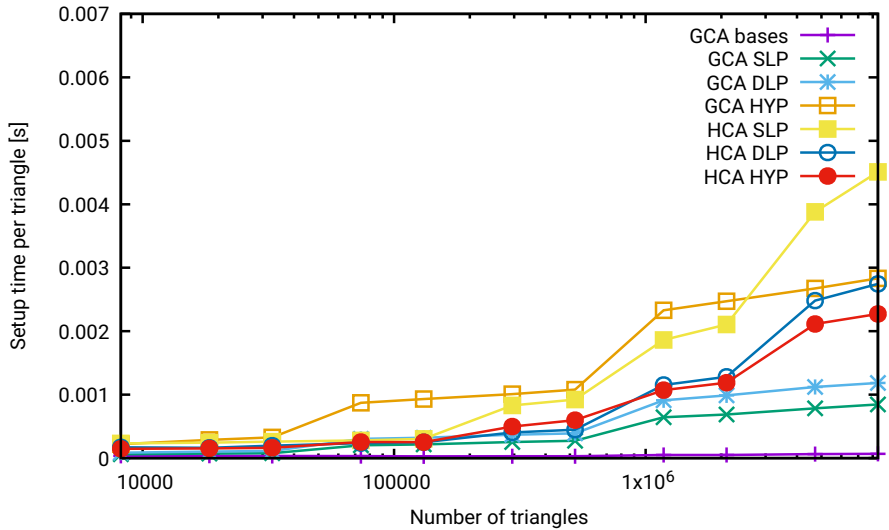

  \pgfdeclareimage[width=12cm]{setup}{gca_hca_cl_setup}

  \begin{center}
    \pgfuseimage{setup}
  \end{center}

  \caption{Setup times for cluster bases and matrices}
  \label{fi:setup}
\end{figure}

Since we expect the runtime to grow like $\mathcal{O}(n \log^\alpha n)$,
cf. \cite{BOGR04} in combination with \cite[Lemma~3.2]{BO07a} for HCA
and \cite{BOCH14} in combination with \cite[Lemma~3.45]{BO10} for GCA,
where $n$ is the number of triangles and $\alpha>0$, we display the
runtime \emph{per triangle} in Figure~\ref{fi:setup}, using a logarithmic
scale for $n$ and a linear scale for the runtime.
We can see that HCA and GCA behave differently if the problem size
grows: the runtime for HCA shows jumps when the order of interpolation
is increased, while the runtime for GCA shows jumps when the order of the
nearfield quadrature is increased.

This observation underlines a fundamental difference between the two
methods: HCA constructs the full coefficient matrix $S$ for every
matrix block, and the matrix $S$ requires $m^6$ coefficients to be
stored and $\mathcal{O}(k m^3)$ operations for the adaptive cross
approximation, where we expect $k\sim m^2$.
GCA, on the other hand, applies cross approximation only per cluster,
not per block, and the recursive structure of the algorithm ensures
that only indices used in a cluster's children are considered in
their parent.
This explains why GCA is less vulnerable to increases in the order
than HCA.

On the other hand, HCA computes the approximation of a submatrix
by evaluating \emph{single} integrals, cf. (\ref{eq:hca_discrete}),
that can be computed in $\mathcal{O}(q_\text{near}^2)$ operations,
while GCA relies on \emph{double} integrals, i.e., the entries of
the matrix $G$, that require $\mathcal{O}(q_\text{near}^4)$ operations.
This explains why HCA is less vulnerable to increases in the
nearfield quadrature order than GCA.

The setup of the matrices dominates the runtime, e.g., for
$8\,388\,688$ triangles, the matrices $G$, $K$, and $W$ take
$7\,096$, $9\,937$, and $23\,752$ seconds to set up with GCA,
respectively, while the preconditioner for $G$ takes only
$4\,898$ seconds for coarsening and $2\,244$ seconds for the
factorization, with $2\,362$ and $1\,470$ seconds for the
preconditioner for $W$.
Solving the linear system with the preconditioned conjugate gradient
method takes around $750$ seconds for the Dirichlet-to-Neumann problems
and around $500$ seconds for the Neumann-to-Dirichlet problems.

\begin{figure}
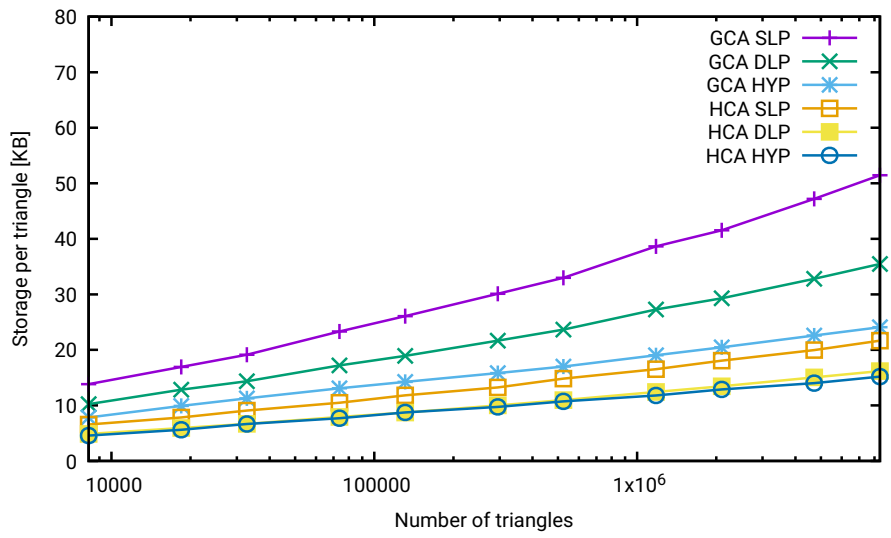

  \pgfdeclareimage[width=12cm]{mem}{gca_hca_cl_mem}

  \begin{center}
    \pgfuseimage{mem}
  \end{center}

  \caption{Storage requirements for the matrices}
  \label{fi:mem}
\end{figure}

Of course, the storage requirements of the algorithms may be even
more important than the runtime, since they determine the size of
a problem that ``fits'' into a given computer.
We again expect a growth like $\mathcal{O}(n \log^\alpha n)$
and report the storage requirements \emph{per triangle}
in Figure~\ref{fi:mem}.

%
%

Although theory leads us to expect the storage requirements
to grow like $\mathcal{O}(n \log^2 n)$, Figure~\ref{fi:mem} suggests
a behaviour more like $\mathcal{O}(n \log n)$ in practice.
We can see that HCA consistently requires less storage than GCA.
This is not surprising, since the algebraic algorithm \cite{BO07a}
employed to turn the hierarchical matrix provided by HCA into an
$\mathcal{H}^2$-matrix essentially computes the \emph{best} possible
$\mathcal{H}^2$-matrix approximation.

We may conclude that a server with 2 processors and just $48=2\times 24$
processor cores equipped with $1\,536$ GB of main memory can handle
boundary element problems with more than 8 million triangles in
a matter of hours without sacrificing accuracy.

\begin{figure}
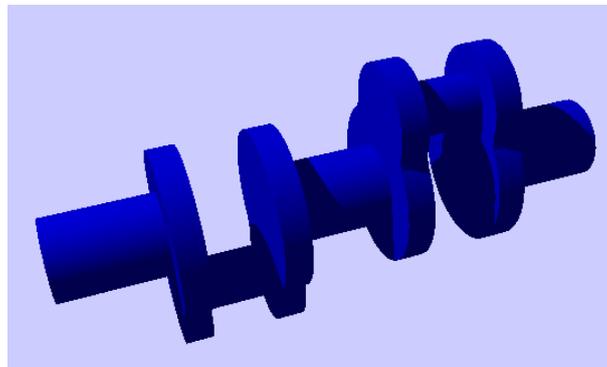

  \pgfdeclareimage[width=8cm]{crankshaft}{sshaft25744}

  \begin{center}
    \pgfuseimage{crankshaft}
  \end{center}

  \caption{Crank shaft geometry}
  \label{fi:crankshaft}
\end{figure}

Admittedly, the unit sphere considered so far is an academic
example.
In order to demonstrate that the techniques also work in more
complicated settings, we consider the crank shaft geometry
displayed in Figure~\ref{fi:crankshaft} created by Joachim Sch\"oberl's
\texttt{netgen} Software.
We start with a mesh with $25\,744$ triangles and refine these triangles
regularly in order to obtain higher resolutions.

%
%

The boundary element mesh is far less ``smooth'' in this case, and
this leads both to an increased condition number and the need to
use significantly higher nearfield quadrature orders.
Since we have already seen that HCA is far less susceptible to
the nearfield quadrature than GCA, we only consider HCA in this
example.
Experiments indicate that the parameters given in Table~\ref{ta:shaft}
are sufficient to preserve the theoretically predicted convergence
rates of the Galerkin method.

%
%
\begin{table}
  \caption{Parameters chosen for the crank shaft geometry}
  \label{ta:shaft}

  \begin{equation*}
    \begin{array}{r|r|r|r|r|r|r|r}
      n & q_\text{near} & r_\text{leaf} & m
        & \epsilon_\text{aca} & \epsilon_\text{comp}
        & \epsilon_\text{slv} & \epsilon_\text{prc}\\
      \hline
     25\,744 & 7 & 64 & 5 & 1_{-9} & 1_{-9} & 1_{-11} & 1_{-3}\\
    102\,976 & 8 & 64 & 6 & 1_{-10} & 1_{-10} & 1_{-12} & 5_{-4}\\
    231\,696 & 9 & 64 & 7 & 1_{-11} & 1_{-11} & 1_{-13} & 2_{-4}\\
    411\,904 & 9 & 64 & 8 & 1_{-12} & 1_{-12} & 1_{-14} & 1_{-4}\\
    926\,784 & 10 & 64 & 9 & 1_{-12} & 1_{-12} & 1_{-14} & 5_{-5}\\
 1\,647\,616 & 10 & 64 & 10 & 1_{-13} & 1_{-13} & 1_{-15} & 5_{-5}\\
    \end{array}
  \end{equation*}
\end{table}

%
%
\begin{figure}
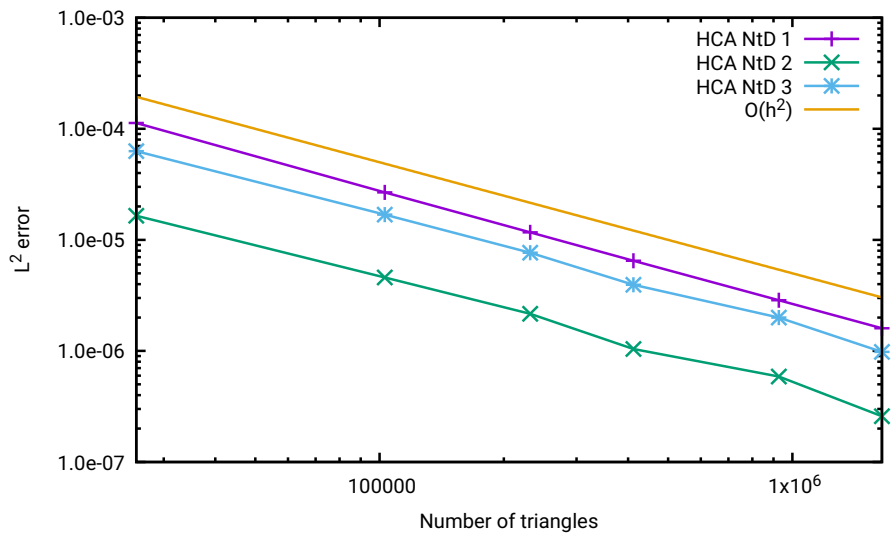

  \pgfdeclareimage[width=12cm]{l2shaft_ntd}{hca_cl_l2shaft_ntd}

  \begin{center}
    \pgfuseimage{l2shaft_ntd}
  \end{center}

  \caption{$L^2$ error for the Neumann-to-Dirichlet problem for the
           crank shaft geometry}
  \label{fi:l2shaft_ntd}
\end{figure}

Figure~\ref{fi:l2shaft_ntd} shows the $L^2$-norm errors for the
Neumann-to-Dirichlet problem at different refinement levels of the
crank shaft geometry.
We can see that the optimal $\mathcal{O}(h^2)$ rate of convergence is
again preserved despite the matrix compression.

%
%
\begin{figure}
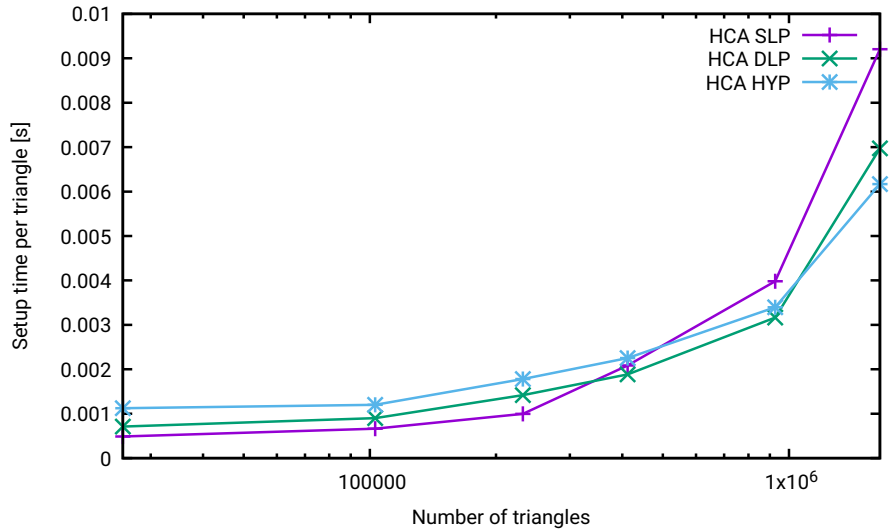

  \pgfdeclareimage[width=12cm]{setupshaft}{hca_cl_setupshaft}

  \begin{center}
    \pgfuseimage{setupshaft}
  \end{center}

  \caption{Setup times for matrices for the crank shaft geometry}
  \label{fi:setupshaft}
\end{figure}

Figure~\ref{fi:setupshaft} shows the setup times per triangle
for the three matrices.
Due to the computationally expensive nearfield quadrature, the
setup time for the matrices dominates the other parts of the
program, e.g., for $1\,647\,616$ triangles the setup times for
the three matrices are $15\,167$, $11\,488$, and $10\,164$
seconds, respectively, while computing both preconditioners takes
only $3\,156$ seconds and each linear system is solved in
under $180$ seconds.

We conclude that using modern compression techniques like
HCA and GCA in combination with efficient $\mathcal{H}^2$-matrix
representations of the resulting matrices, large boundary
element problems on meshes with several million triangles
can be treated in few hours on moderately expensive servers.

\bibliographystyle{splncs04}
\bibliography{hmatrix}

\end{document}